\documentclass[11pt]{article}
\usepackage{amssymb}
\usepackage{latexsym}
\usepackage{amsmath}
\usepackage{amsthm}
\usepackage{mathtools}
\usepackage[round,sort]{natbib}
\usepackage{tikz-cd}
\usepackage{enumitem} 
\usepackage{hyperref}
\usepackage{abstract}
\hypersetup{
    colorlinks,
    citecolor=black,
    filecolor=black,
    linkcolor=black,
    urlcolor=black
}
\newtheorem{thm}{Theorem}[section]

\newtheorem{lemma}[thm]{Lemma}
\newtheorem{cor}[thm]{Corollary}
\newtheorem{dfn}[thm]{Definition}

\newcommand{\reals}{\mathbb R}
\newcommand{\cplx}{\mathbb C}
\newcommand{\intg}{\mathbb Z}
\newcommand{\ratl}{\mathbb Q}

\newcommand{\calo}{{\cal O}}

\renewcommand{\qedsymbol}{$\square$}
\title{\bf Witten Genus on String Toric Complete Intersections}
\author{Lin-Da Xiao}
\date{\vspace{-5ex}}
\begin{document}
\maketitle
\begin{abstract}
By using the equivariant localization formula of toric varieties. We prove the vanishing of the Witten genus of some string complete intersections in smooth toric varieties.
\end{abstract}
\section{Introduction}
Let $M$ be a $4k$ dimensional closed oriented smooth manifold. In \citet{witten1988index}, a multiplicative genus is defined by formally applying the equivariant Atiyah-Singer index theorem to a hypothetical Dirac operator in the loop space and we can get the analogue of $\hat{A}$-genus
\begin{equation*}
	W(M)=\left\langle\hat{A}(TM) Ch(\Theta(T_\mathbb{C}M)),[M] \right\rangle,
\end{equation*}
where 
\begin{equation*}
	\Theta(T_{\mathbb{C}}M)=\underset{m=1}{\overset{\infty}{\bigotimes}} S_{q^{m}}(T_{\mathbb{C}}M-\mathbb{C}^{4k})
\end{equation*}
is the Witten bundle where $q=e^{2\pi\sqrt{-1} \tau}$ with $Im(\tau)\geq 0$. Also, 
		\begin{equation*}
			S_{q^m} (T_{\mathbb{C}}M-\mathbb{C}^{4k}) :=\sum^\infty_{k=0}(S^k (T_{\mathbb{C}}M-\mathbb{C}^{4k}))(q^m)^k,
		\end{equation*} 
where $S^k(T_\cplx M-\cplx^{4k})$is the k-th symmetric power of the formal difference $T_\cplx M-\cplx^{4k}$ in the K-theory.
To manifest the modular aspects of Witten genus, following \cite{liu1996elliptic}, we can also write it as
\begin{equation*}
	W(M)=\left\langle\prod_{i}\frac{z_i \theta'(0,\tau)}{\theta(z_i,\tau)},[M]\right\rangle
\end{equation*}
where $\{\pm 2\pi\sqrt{-1}z_i,1\leq i\geq 2k\}$ are the formal Chern roots of the bundle $T_{\mathbb{C}}M$.

The oriented manifold $M$ is called \textbf{spin} if the  second Stiefel-Whitney class $w_2(TM)$ vanishes. Moreover, manifold $M$ is called \textbf{string} if the half first Pontryagin class vanishes. According to the Atiyah-Singer index theorem, when manifold is spin, the Witten genus is an integral expansion in terms of $q$ (cf. \cite{hirzebruch1992manifolds}). It is also well known that if the manifold is string, the Witten genus is a modular form of weight $2k$ over $SL(2,\mathbb{Z})$ with integral Fourier expansion (\cite{zagier1988note}).
Analogous to Lichnerowicz's classical result on $\hat{A}$ genus (cf. \cite{lawson2016spin}), which stated that $\hat{A}$ genus on spin manifold with positive scalar curvature vanishes, Stolz conjectured that the Witten genus on string manifold with positive Ricci curvature vanishes(for the original arguments, cf. \cite{stolz1996conjecture}, and for a review cf. \cite{dessai2009some}). 

Toward the Stolz's conjecture, several vanishing results have been discovered. There are basically two types of methods. One is to apply the theorem of \cite{dessai1994witten} which is based on \cite{liu1995modular} when the manifold admits some nontrivial action of a semi-simple Lie group. The current results via this method include:
\begin{enumerate}
\item String homogeneous spaces of compact semi-simple Lie groups.(\cite{dessai1994witten,liu1995modular}.
\item Total spaces of fiber bundles with fiber $G/H$, with compact semi-simple structure group $G$ \cite{stolz1996conjecture}.
\item Generalized string complete intersection in irreducible, compact, Hermitian, symmetric spaces\cite{forster2007stolz}.
\item String manifold with effective torus action such that $dim\  T>b_2(M)$ \cite{wiemeler2017note}, where $T$ is a compact torus and $b_2(M)$ is the second Betti number of $M$.
\end{enumerate}

Alternatively, sometimes, one can reduce the calculation of Witten genus to calculation of residues. This method was first used by Landweber and Stone \cite{hirzebruch1992manifolds}, which is purely computational and more direct. The vanishing results include:
\begin{enumerate}
\item String complete intersection in projective space\cite{hirzebruch1992manifolds}.
\item String generalized complete intersection in products of projective spaces \cite{chen2008witten}。
\item String complete intersection in products of Grassmannians and flag manifolds\cite{zhou2014witten,zhuang2016vanishing}。
\end{enumerate}
This method also has many applications in elliptic genus, e.g., \cite{ma2005elliptic,gorbounov2008mirror}.
In this paper we generalize the result of~\cite{chen2008witten} to string complete intersections in Toric varieties. We mainly follow the second method in our calculation and also borrow the techniques of equivariant localization in \cite{dessai2016torus}.

\subsection*{Main Result}
Consider a symplectic toric variety $X$ with a set of  invariant divisors $\{D_{\rho_j}\}_{1\leq j\leq r}$ and of Picard number $k$. Choose a basis $\{q_1,...,q_k\}$ of the Picard group to expand all the invariant divisors with integer coefficients
\begin{equation*}
D_{\rho_j}=\left\{
\begin{aligned}
& q_j & 1\leq i\leq k\\
& \sum_{i=1}^k m_{j i} q_i  &k+1\leq j\leq r.
\end{aligned}
\right.
\end{equation*}
Consider a smooth complete intersection $Y\subset X$ as a generic intersection of $s$ hypersurfaces $\{Y_l\}_{1\leq l\leq s}$. Each $Y_l$ is dual to the cohomology class $\sum_{j}^k d_{l i} q_i$.
\begin{thm}\label{thm:main}
When the integer matrix elements $(m_{ji})$ and $(d_{ji})$ satisfy
$$
\sum_{j=1}^n d_{ji} d_{j l}-\sum_{j=k+1}^r m_{ji}m_{j l}=0 \text{  for } i\neq l.
$$
and
$$
\sum_{j=1}^n d_{ji}^2-\sum_{j=k+1}^r m_{ji}^2-1=0,
$$
then the complete intersection $Y$ is string and its Witten genus vanishes.
\end{thm}
\section{Preparations}
We will gather all the necessary notions and lemmas in this section.
\begin{dfn}
$M$ is a $4k$-dimensional compact oriented manifold with a real vector bundle $E$  of rank $2n$, the \textbf{Witten class} of $E$ is defined to be 
$$
\mathcal{W}(E,M):=\hat{A}(E)Ch(\Theta(E\otimes \cplx)).
$$ 
\end{dfn}
One can check that Whitney product formula holds for Witten class, i.e. for an exact sequence of real vector bundle $0\rightarrow E\rightarrow F\rightarrow G\rightarrow 0$, we have
$$
\mathcal{W}(E)\cdot\mathcal{W}(G)=\mathcal{W}(F).
$$

It is customary to choose Chern roots $\{\pm 2\pi\sqrt{-1}z_i,1\leq i\geq 2k\}$, which simplifies the Witten genus to 
\begin{equation*}
	W(M)=\left\langle\prod_{i}\frac{z_i \theta'(0,\tau)}{\theta(z_i,\tau)},[M]\right\rangle,
\end{equation*}
where $\theta(x,\tau)$ is the first Jacobi theta function.
\subsection*{Global Residue Theorem}
First, we recall the residue theorem (Chapter 5 of~\cite{Griffiths1994principles}). Let $M$ be a compact complex manifold of dimension $s$. Suppose that $D_i$ for $i=1,...,m$ are effective divisors, the intersection of which is a finite set of points. Let $D=D_1+...+D_m$. Let $\omega$ be a meromorphic $m-$form on $M$ with polar divisor $D$. For each point $P\in D_1\cap...\cap D_m$, we may restrict $\omega$ to a neighborhood $U_P$ of $P$ and define the residue at $P$, denoted by $Res_P\omega$. Then, one has:
\begin{lemma}(Residue theorem).  
$$
\sum_{P\in D_1\cap...\cap D_m} Res_P\omega=0.
$$
\end{lemma}
\subsubsection*{Jacobi Theta Functions}
We collect all the necessary fact of Jacobi theta function here, all of which can be found in \cite{chandrasekharan1985elliptic}. The theta functions are defined as follows:
\begin{equation*}
\begin{aligned}
	\theta (v,\tau ) & =2 q^{1/8} \sin(\pi  v) \prod _{j=1}^{\infty } \left[\left(1-q^j\right) \left(1-e^{2 \pi  i v} q^j\right) \left(1-e^{-2 \pi  i v} q^j\right)\right],\\
	\theta _1(v,\tau )& =2 q^{1/8} \cos(\pi  v) \prod _{j=1}^{\infty } \left[\left(1-q^j\right) \left(1+e^{2 \pi  i v} q^j\right) \left(1+e^{-2 \pi  i v} q^j\right)\right],\\
\end{aligned}
\end{equation*}
\begin{equation*}
\begin{aligned}
	\theta _2(v,\tau )&=\prod _{j=1}^{\infty } \left[\left(1-q^j\right) \left(1-e^{2 \pi  i v} q^{j-\frac{1}{2}}\right) \left(1-e^{-2 \pi  i v} q^{j-\frac{1}{2}}\right)\right],\\
	\theta _3(v,\tau )&=\prod _{j=1}^{\infty } \left[\left(1-q^j\right) \left(1+e^{2 \pi  i v} q^{j-\frac{1}{2}}\right) \left(1+e^{-2 \pi  i v} q^{j-\frac{1}{2}}\right)\right],
\end{aligned}
\end{equation*}
where $q=e^{2\pi\sqrt{-1} \tau}$ with $Im(\tau)\geq 0$. We also have the Jacobi identity:
\begin{equation*}
\theta '(0,\tau )=\frac{\partial }{\partial v}\theta (v,\tau )|_{v=0}=\pi\theta_1(0,\tau )\theta _2(0,\tau )\theta _3(0,\tau ).
\end{equation*}
They satisfy the transformation law under the translation on lattice $\{\intg+b\intg\}$
\begin{equation}\label{theta}
\begin{aligned}
\ &\theta (v+m,\tau ) =(-1)^m\theta (v,\tau ),& &\theta (v+n \tau ,\tau )=(-1)^n e^{ -2 \pi  i n v-\pi  i n^2 \tau}\theta (v,\tau ),\\
\ &\theta_1 (v+m,\tau ) =(-1)^m\theta_1 (v,\tau ),& &\theta_1 (v+n \tau ,\tau )= e^{ -2 \pi  i n v-\pi  i n^2 \tau}\theta_1 (v,\tau ),\\
\ &\theta_2 (v+m,\tau ) =\theta_2 (v,\tau ), & &\theta_2 (v+n \tau ,\tau )=(-1)^n e^{ -2 \pi  i n v-\pi  i n^2 \tau}\theta_2 (v,\tau ),\\
\ &\theta_3 (v+m,\tau ) =\theta_3 (v,\tau ),& &\theta_3 (v+n \tau ,\tau )= e^{ -2 \pi  i n v-\pi  i n^2 \tau}\theta_3 (v,\tau ).
\end{aligned}
\end{equation}
\subsection*{Equivariant Localization and Genera as Residues}

There is a ring morphism between the ordinary cohomology ring $H^\bullet(X,\mathbb{Q})$ and the equivariant cohomology ring $H_{T}^{\bullet}(X,\mathbb{Q})$. $H^\bullet_T(X_\Sigma)$ can be thought of as an $H^\bullet(BT)$-module. Note also that the inclusion of a fiber $i_{X}: X_\Sigma\rightarrow X\times_G EG$ induces the ``Non-equivariant limit'' map $i^*_X: H^\bullet_T(X_\Sigma)\rightarrow H^\bullet(X_\Sigma)$, which amounts to mapping all $\lambda_i$ to 0. In the case of toric varieties, there is a equivariant version of Jurkiewicz-Danilov Theorem, see Chapter 12 of \cite{cox2009toric}
\begin{lemma}
Every $D_\rho$ has an equivariant counterpart $(D_\rho)_T$, where $\rho\in \Sigma(1)$. This map induces the morphism between cohomology and their equivariant counterpart.
\begin{equation*}
\begin{aligned}
D_\rho & \longrightarrow (D_{\rho})_T=D_\rho-\lambda_\rho\in H^2_T(X,\ratl)\\
 H^\bullet (X,\ratl)&\longrightarrow H^\bullet_T(X_\Sigma,\mathbb{Q})=\mathbb{Q}[(D_{\rho_1})_T,...,(D_{\rho_r})_T]/(\mathcal{I}_T+\mathcal{J}_T)\\
\end{aligned}
\end{equation*}
where
\begin{equation*}
\mathcal{I}_T = \langle (D_{i_1})_T\cdot\cdot\cdot (D_{i_s})_T | i_j \text{ are distinct and } \rho_{i_1}+...+\rho_{i_s} \text{ is not a cone of } \Sigma\rangle
\end{equation*}
and $\mathcal{J}_T$ is the ideal generated by the linear forms
\begin{equation*}
\sum_{i}^r\langle m, u_i \rangle (D_{\rho_i})_T
\end{equation*}
where $m$ ranges over $M$ (or equivalently, over some basis for $M$).
\end{lemma}

We can find an explicit way to proceed calculations via the Atiyah-Bott localization. The following formula can be found in \cite{givental1998mirror} 
\begin{lemma}
$f(q_1,...,q_k,\{\lambda_i\})$ is a polynomial in $H^\bullet_T(X,\ratl)$, integrating it over the fundamental class $[X]$ maps it into $\ratl[\lambda_1,...,\lambda_r]$. Explicitly, we have
\begin{equation}\label{localization}
\int_{X} f(q_1,...,q_k,\{\lambda_i\})=\sum_\alpha Res_\alpha \frac{f(q_1,...,q_k,\{\lambda_i\})}{(D_1)_T\cdot (D_2)_T\cdot ...\cdot (D_r)_T} d q_1 dq_2...dq_k.
\end{equation}
The symbol $Res_\alpha$ refers to the residue of the $k-$form at the pole by the order set of equations
$$
\sum^k_{i=1}p_i m_{i j_s}=\lambda_s,\ s=1,...,k 
$$
\end{lemma}

\section{Main Results}
From now on, we will abuse the notation to write $D_\rho$ instead of $[D_\rho]$ as the cohomology classes and divisor classes. We consider a smooth compact toric variety $X_\Sigma$ corresponding to $\Sigma\in N_\reals$, where $N$ is a lattice of dimension $n$ and $M$ being its dual lattice. Let $\Sigma(1)$ denote the set of one dimensional cones in $\Sigma$ and assume that $|\Sigma(1)|=r$. By Proposition 4.2.1 and 4.2.5 in \cite{cox2009toric}, we know $\text{Pic}(X_\Sigma)$ is torsion-free and
\begin{equation*}
0\longrightarrow M\longrightarrow \intg^r\longrightarrow \text{Pic}(X_\Sigma)\longrightarrow 0.
\end{equation*}
If we set the Picard number $k$, we have the relation
\begin{equation*}
	\text{dim }X_\Sigma =\text{dim }N=n=r-k.
\end{equation*}
After we quotient the ideal generated by the linear relations, we can equivalently say that the cohomology ring is multiplicatively generated by a basis of the Picard group. We denote the basis of Picard group as $q_1,...q_k$, and express the invariant divisors as
\begin{equation*}
D_{\rho_i}=\left\{
\begin{aligned}
& q_i & 1\leq i\leq k\\
& \sum_{j=1}^k m_{i j} q_j  &k+1\leq i\leq r
\end{aligned}
\right.
\end{equation*}
Then let's consider a smooth complete intersection $Y$ as a generic intersection of $s$ hypersurfaces $\{Y_l\}_{1\leq l\leq s}$. Each $Y_l$ is dual to the cohomology class $\sum_{j}^k d_{l j} q_j$. In the context of smooth manifold, by the transversality argument, one can always deform the hypersurfaces (not algebraic now) so that they intersects transversally. Thus $Y$ is submanifold of complex dimension $r-k-s$. In this section we will consider the Witten genus on $Y$. From now on, we will use both $X:=X_\Sigma$ to denote the toric variety, when there is no confusion.
\subsection{Proof of Theorem~\ref{thm:main}}
For smooth toric varieties, there is a \textbf{generalized Euler sequence}( see Theorem 8.1.6 in \cite{cox2009toric}).
\begin{equation*}
0\longrightarrow \mathcal{O}_{X_\Sigma}^{\oplus k}\longrightarrow\oplus_{\rho\in \Sigma(1)}\mathcal{O}_{X_\Sigma}(D_\rho)\longrightarrow TX_{\Sigma}\longrightarrow0,
\end{equation*}
where $k$ is the rank of $\text{Pic}(X_\Sigma)$ or equivalently the rank of $H^2(X_\Sigma,\mathbb{Z})$.
from which we can calculate the total Chern class and all the multiplicative genera of $TX_{\Sigma}$. 

The $\mathcal{O}_{X_\Sigma}(D_\rho)$ is the line bundle defined by the divisor $D_\rho$ in the Picard group, of which the first Chern class is $D_\rho\in H^2(X_\Sigma,\intg)$.
Then by the multiplicative property of total Chern class, we have
\begin{equation*}
	c(TX_\Sigma)=\oplus_{\rho\in\Sigma(1)}\  c(\mathcal{O}_{X_\Sigma}(D_\rho))=\prod_{\rho\in \Sigma(1)}(1+D_\rho).
\end{equation*}
The $\{D_\rho\}$ play the role of Chern roots here, which is not coincidence.
For a general multiplicative genus $\varphi_Q$ corresponding to even power series $Q(x)$, we observe that
\begin{equation*}
\varphi_Q(TX_\Sigma)=\frac{\varphi_Q\left(\bigoplus_{\rho\in \Sigma(1)}\mathcal{O}(D_\rho)\right)}{\varphi_Q(\mathcal{O}(X_\Sigma)^{\oplus k})}=\frac{\prod_{\rho\in \Sigma(1)} Q(D_\rho)}{Q(0)^k},
\end{equation*}
In the case of Witten genus, we have
\begin{equation*}
\begin{aligned}
\mathcal{W}(T_\reals X)&=\prod_{\rho\in \Sigma(1)}D_\rho \frac{\theta'(0,\tau)}{\theta(D_\rho,\tau)},\\
W(X_\Sigma)&=\left\langle\mathcal{W}(X_\Sigma),[X_\Sigma]\right\rangle
\end{aligned}
\end{equation*}
Consider the inclusion map $\iota: Y\rightarrow X_\Sigma$, we have the adjunction formula:
$$
0\longrightarrow TY\longrightarrow\iota^*TX_{\Sigma}\longrightarrow\iota^* N_Y X\longrightarrow0,
$$
where the normal bundle $N_Y X$ is isomorphic to $\bigoplus_{l=1}^s\mathcal{O}(\sum_{j}^k d_{l j} q_j)$.
By the multiplicative properties of Chern class, we have the following data:
\begin{equation*}
c(TY)=\iota^*\left(\frac{c(TX)}{c(N_Y X)}\right)=\iota^*\left(\frac{\prod_{i=1}^k (1+q_i)\prod_{j=k+1}^{r}(1+\sum_{i=1}^k m_{ji}q_i)}{\prod_{l=1}^s(1+\sum_{j}^k d_{l j} q_j)}\right).
\end{equation*}
Following these we have
$$
\begin{aligned}
w_2(T_\mathbb{R}Y)& \equiv c_1(TY) (\text{mod } 2)\\
& =\iota^*\left(\sum_{i=1}^k q_i+\sum_{j=k+1}^{r} \sum_{i=1}^k m_{ji}q_i-\sum_{l=1}^s\sum_{i=1}^k d_{l i}q_i\right) \text{     (mod } 2)
\end{aligned}
$$
and 
$$
p_1(T_{\mathbb{R}}Y)=\iota^*\left(\sum_{i=1}^k q_i^2+\sum_{j=k+1}^r (\sum_{i=1}^k m_{ji}q_i)^2-\sum_{l=1}^s(\sum_{i=1}^k d_{l i}q_i)^2\right).
$$
Because we know little about the $H^4(X,\intg)$, it is difficult to tell whether the above presentation of Pontryagin class vanishes in the general case. Still, we know the complete intersection is string when all the coefficients of $\{q_i q_j\}$ vanishes. On the other hand, because of the multiplicative property of Witten class, we have
\begin{equation*}
\begin{aligned}
\mathcal{W}(T_\reals Y)
&=\iota^*\left(\frac{\mathcal{W}(T_\reals X)}{\mathcal{W}((N_Y X)_\reals)}\right)\\
&=\iota^*\left(\frac{\prod_{i=1}^k q_i \frac{\theta'(0,\tau)}{\theta(q_i,\tau)}\prod_{j=k+1}^r (\sum_{t=1}^k m_{j t}q_t) \frac{\theta'(0,\tau)}{\theta(\sum_{t=1}^k m_{jt}q_t,\tau)}}{\prod_{l=1}^s(\sum_{j=1}^k d_{l j} q_j) \frac{\theta'(0,\tau)}{\theta(\sum_{j=1}^k d_{l j} q_j,\tau)}}\right).
\end{aligned}
\end{equation*}
Integrate the above formula over the fundamental class of Y, we have
$$
\begin{aligned}
&\int_{Y} \mathcal{W}(T_\reals Y) \\
& =\int_X \iota_!\iota^*\left(\frac{\mathcal{W}(T_\reals X)}{\mathcal{W}((N_Y X)_\reals)}\right)\\
& =\int_X \left(\frac{\mathcal{W}(T_\reals X)}{\mathcal{W}((N_Y X)_\reals))}\right)\smile e(N_Y X)\\
& =\int_X \left(\frac{\prod_{i=1}^k q_i \frac{\theta(q_i,\tau)}{\theta'(0,\tau)}\prod_{j=k+1}^r (\sum_{t=1}^k m_{j t}q_t) \frac{\theta(\sum_{t=1}^k m_{jt}q_t,\tau)}{\theta'(0,\tau)}}{\prod_{l=1}^s(\sum_{j=1}^k d_{l j} q_j) \frac{\theta(\sum_{j=1}^k d_{l j} q_j,\tau)}{\theta'(0,\tau)}}\right)\cdot \prod_{l=1}^s(\sum_{j=1}^k d_{l j} q_j)\\
& =\int_X \left(\frac{\prod_{i=1}^k q_i \frac{\theta(q_i,\tau)}{\theta'(0,\tau)}\prod_{j=k+1}^r (\sum_{t=1}^k m_{jt}q_t) \frac{\theta(\sum_{t=1}^k m_{jt}q_t,\tau)}{\theta'(0,\tau)}}{\prod_{l=1}^s \frac{\theta(\sum_{j=1}^k d_{l j} q_j,\tau)}{\theta'(0,\tau)}}\right),
\end{aligned}
$$
where $e(N_Y X)$ in the second line is the Euler class of the normal bundle.

If we take the non-equivariant limit we can get an explicit formula to calculate the the integration of $f(p,0)$ over the fundamental class $[X]$.
\subsection{Vanishing result for some string complete intersection}
What we want to calculate is
\begin{equation*}
\begin{aligned}
\int_{Y} \mathcal{W}(T_\reals Y) 
& =\int_X \iota_!\iota^*\left(\frac{\mathcal{W}(T_\reals X)}{\mathcal{W}((N_Y X)_\reals)}\right)\\
& =\int_X \left(\frac{\prod_{i=1}^k q_i \frac{\theta'(0,\tau)}{\theta(q_i,\tau)}\prod_{j=k+1}^r (\sum_{t=1}^k m_{jt}q_t) \frac{\theta'(0,\tau)}{\theta(\sum_{t=1}^k m_{jt}q_t,\tau)}}{\prod_{l=1}^s \frac{\theta'(0,\tau)}{\theta(\sum_{j=1}^k d_{l j} q_j,\tau)}}\right).
\end{aligned}
\end{equation*}
Note that the complete intersection $Y$ is generally not invariant under the action of torus, but remember that $TX_\Sigma$ and $N_Y X\cong \calo(\sum^k_j d_{lj} q_j)$ are invariant bundles, we know that the differential form on the right hand side is invariant.
With the localization argument before, we can consider a integral of polynomial $\mathcal{W}(q,\{\lambda_i\})$ in $H^\bullet_T(X,\ratl)$ which under nonequivariant limit gives the integral we want. The choice is not unique. We can choose

\begin{equation*}
\begin{aligned}
&\int_{Y} \mathcal{W}(T_\reals Y)\\ 
& =i_X^*\int_X \left(\frac{\prod_{i=1}^k (q_i-\lambda_i) \frac{\theta'(0,\tau)}{\theta(q_i-\lambda_i,\tau)}\prod_{j=k+1}^r (\sum_{t=1}^k m_{jt}q_t-\lambda_j) \frac{\theta'(0,\tau)}{\theta(\sum_{t=1}^k m_{jt}q_t-\lambda_j,\tau)}}{\prod_{l=1}^s \frac{\theta'(0,\tau)}{\theta(\sum_{j=1}^k d_{l j} q_j,\tau)}}\right)\\
\end{aligned}
\end{equation*}
Then by Eq.~\ref{localization},
\begin{equation*}
\begin{aligned}
&\int_{Y} \mathcal{W}(T_\reals Y) \\
& =i_X^*\sum_\alpha res_\alpha \frac{\left(\frac{\prod_{i=1}^k (q_i-\lambda_i) \frac{\theta'(0,\tau)}{\theta(q_i-\lambda_i,\tau)}\prod_{j=k+1}^r (\sum_{t=1}^k m_{jt}q_t-\lambda_j) \frac{\theta'(0,\tau)}{\theta(\sum_{t=1}^k m_{jt}q_t-\lambda_j,\tau)}}{\prod_{l=1}^s \frac{\theta'(0,\tau)}{\theta(\sum_{j=1}^k d_{l j} q_j,\tau)}}\right)}{\prod_{i=1}^k (q_i-\lambda_i)\prod_{j=k+1}^r (\sum_{t=1}^k m_{jt}q_t-\lambda_j)}\\
& =i_X^*\sum_\alpha res_\alpha\left(\frac{\prod_{l=1}^s \frac{\theta(\sum_{j=1}^k d_{l j} q_j,\tau)}{\theta'(0,\tau)}}{\prod_{i=1}^k  \frac{\theta(q_i-\lambda_i,\tau)}{\theta'(0,\tau)}\prod_{j=k+1}^r  \frac{\theta(\sum_{t=1}^k m_{jt}q_t-\lambda_j,\tau)}{\theta'(0,\tau)}}\right)dq_1...dq_k\\
&=\lim_{\lambda_i\rightarrow0}\sum_\alpha res_\alpha \frac{g(q,\tau)}{\prod_{j=1}^r f_j (q,\tau,\{\lambda_i\})}dq_1...dq_k,
\end{aligned}
\end{equation*}
where
\begin{equation*}
\begin{aligned}
g(q_1,..,q_k,\tau)=\prod_{l=1}^s\frac{\theta(\sum_{j=1}^k d_{l j} q_j,\tau)}{\theta'(0,\tau)}
\end{aligned}
\end{equation*}
and
\begin{equation*}
f_j(q_1,...,q_k,\tau,\{\lambda_i\})=\left\{
\begin{aligned}
&\frac{\theta(q_j-\lambda_j,\tau)}{\theta'(0,\tau)} & \text{ for } 1\leq j\leq k\\
& \frac{\theta(\sum_{t=1}^k m_{jt}q_t-\lambda_j,\tau)}{\theta'(0,\tau)} &\text{for } k+1\leq j\leq r
\end{aligned}
 \right.
\end{equation*}
Let's analyze how $g$ and $f$ transform under the translation over the lattice $\{\intg+\intg \tau\}$. By the transforming laws of theta functions[\ref{theta}], without loss of generality, $q_1\rightarrow q_1+1$.
\begin{equation*}
	g(q_1+1,...,q_k,\tau)=(-1)^{d_{11}+...+d_{s1}}g(q_1,...,q_k,\tau),
\end{equation*}
\begin{equation*}
f_j(q_1+1,...,q_k,\tau,\{\lambda_i\})=\left\{
\begin{aligned}
	& - f_1(q_1,...,q_k,\tau,\{\lambda_i\})&\text{for $j=1$}\\
	& f_j(q_1,...,q_k,\tau,\{\lambda_i\})&\text{for $j \neq 1$ and $j\leq k$} \\
	& (-1)^{m_{j1}} f_j (q_1,...,q_k,\tau,\{\lambda_i\})&\text{for $k+1\leq j\leq r$}
\end{aligned}
\right.
\end{equation*}
Then,
\begin{equation*}
\begin{aligned}
&\frac{g(q_1,...,q_i+1,...,q_k,\tau)}{\prod_{j=1}^r f_j (q_1,...,q_i+1,...,q_k,\tau,\{\lambda_i\})}\\
&=(-1)^{\sum_{j=1}^s d_{ji}-\sum_{j=k+1}^r m_{ji}-1}\frac{g(q_1,...,q_k,\tau)}{\prod_{j=1}^r f_j (q_1,...,q_k,\tau,\{\lambda_i\})}.
\end{aligned}
\end{equation*}
On the other hand under the translation $q_1\rightarrow q_1+\tau$, we have
\begin{equation*}
g(q_1+\tau,...,q_k,\tau)
=(-1)^{d_{11}+...,d_{s1}} e^{-2\pi i\sum_{l=1}^s d_{l1}(\sum^k_j d_{lj}q_j)-\pi i\sum_{l=1}^s d^2_{l1}\tau} g(q_1,...,q_k,\tau),
\end{equation*}
and
\begin{equation*}
\begin{aligned}
&f_j(q_1+\tau,...,q_k,\tau,\{\lambda_i\})=\\
&\left\{
\begin{aligned}
&(-1)e^{-2\pi i (q_1-\lambda_1)-\pi i \tau} f_1(q,\tau,\{\lambda_i\}) &\text{ for $j=1$ }\\
&f_j(q,\tau,\{\lambda_i\}) &\text{ for $i\neq 1$ and $j\leq k$}\\
& (-1)^{m_{j1}} e^{-2\pi i m_{j1}(\sum^k_t m_{jt}q_t-\lambda_j)-\pi i m_{j1}^2 \tau}f_j(q,\tau,\{\lambda_i\})&\text{$k+1 \leq j\leq r$}
\end{aligned}
\right.
\end{aligned}
\end{equation*}
Then the we have the transforming law
\begin{equation*}
\begin{aligned}
&\frac{g(q_1,...,q_i+\tau,...,q_k,\tau)}{\prod_{j=1}^r f_j (q_1,...,q_i+\tau...,q_k,\tau,\{\lambda_i\})}\\
&=(-1)^{\sum_{t=1}^s d_{ti}-\sum_{t=k+1}^r m_{ti}-1}\\
&\times e^{-2\pi \sqrt{-1}[\sum_{v\neq i}^k(\sum_{t}^s d_{t i} d_{tv}-\sum_{u=k+1}^r m_{u i} m_{uv})q_v+(\sum_{t=1}^s d_{ti}^2-\sum^r_{u=k+1}m_{ui}^2-1)q_i-\lambda_i+\sum_{u=k+1}^r m_{ui}\lambda_u]}\\
&\times e^{-\pi \sqrt{-1}(\sum_{t=1}^s d_{ti}^2-\sum^r_{u=k+1}m_{ui}^2-1)\tau}\\
&\times \frac{g(q_1,...,q_i,...,q_k,\tau)}{\prod_{j=1}^r f_j (q_1,...,q_i...,q_k,\tau,\{\lambda_i\})}
\end{aligned}.
\end{equation*}
Remember that
\begin{equation*}
\begin{aligned}
w_2(T_\mathbb{R}Y)& \equiv c_1(TY) (\text{mod } 2)\\
& =\iota^*\left(\sum_{i=1}^k q_i(1+\sum_{j=k+1}^{r} m_{ji}-\sum_{l}^s d_{l i})\right) \text{     (mod } 2),
\end{aligned}
\end{equation*}
and
\begin{equation*}
\begin{aligned}
p_1(T_{\mathbb{R}}Y)&=\iota^*\left(\sum_{i=1}^k q_i^2+\sum_{j=k+1}^r (\sum_{i=1}^k m_{ji}q_i)^2-\sum_{l=1}^s(\sum_{i=1}^k d_{l i}q_i)^2\right)\\
&=\iota^*\left(\sum_{i=1}^k q_i^2(1+\sum_{j=k+1}^r m_{ji}^2-\sum_{l=1}^s d_{l i}^2)+\sum_{i=1}^k\sum_{l=1}^k q_i q_l(\sum^r_{j=k+1}m_{ji}m_{jl}-\sum_{u}^s d_{u i}d_{ul})\right)\\
\end{aligned}
\end{equation*}
Also notice that, when
\begin{equation*}
1+\sum_{j=k+1}^r m_{ji}^2-\sum_{l=1}^s d_{l i}^2=0,
\end{equation*}
$\omega_2(T_\reals Y)$ vanishes automatically because
\begin{equation*}
\sum_{j=k+1}^r m_{ji}^2-\sum_{l=1}^s d_{l i}^2\equiv \sum_{j=k+1}^{r} m_{ji}-\sum_{l}^s d_{l i}\ (\mod 2).
\end{equation*}
At nonequivariant limit, we take $\lambda_i=0$, assuming
$$
\sum_{j=1}^n d_{ji} d_{j l}-\sum_{j=k+1}^r m_{ji}m_{j l}=0 \text{  for } i\neq l.
$$
and
$$
\sum_{j=1}^n d_{ji}^2-\sum_{j=k+1}^r m_{ji}^2-1=0,
$$
then the complete intersection $Y$ is string and
the function 
\begin{equation*}
\frac{g(q,\tau)}{\prod_{j=1}^r f_j (q,\tau,0)}
\end{equation*}
is elliptic over the lattice $\Gamma:=\{\intg+\intg\tau\}$. Then equivalently, we 
can regard 
\begin{equation*}
\frac{g(q,\tau)}{\prod_{j=1}^r f_j (q,\tau,0)}dq_1...dq_k
\end{equation*}
as a meromorphic form over the compact torus $(\cplx/\Gamma)^s$.
Then by the global residue theorem, the Witten genus
\begin{equation*}
\begin{aligned}
&\int_{Y} \mathcal{W}(T_\reals Y) \\
&=\sum_\alpha res_\alpha \frac{g(q,\tau)}{\prod_{j=1}^r f_j (q,\tau,0)}dq_1...dq_k,
\end{aligned}
\end{equation*}
vanishes.
\begin{flushright}
\qedsymbol
\end{flushright}
\begin{cor}
In the special case where the toric variety is  product of projective spaces of the form $\mathbb{P}^{n_1}\times \mathbb{P}^{n_2}\times\cdot \cdot \cdot \times \mathbb{P}^{n_t}$, then
up to reordering, the matrix $m_{ij}$ can be written as
$$
\begin{pmatrix}
\underbrace{1,1,...,1}_{n_{1}} & 0 & \dots  & 0\\
0 & \underbrace{1,1,...,1}_{n_{2}} & \dots  & 0\\
\vdots  & \vdots  & \ddots  & \vdots \\
0 & 0 & \dots  & \underbrace{1,...,1}_{n_{k}}
\end{pmatrix}.
$$
The vanishing condition reduces to 
$$
\sum_{j=1}^n d_{ji} d_{j l}=0 \text{  for } i\neq l.
$$
and
$$
\sum_{j=1}^n d_{ji}^2-n_i-1=0,
$$
which reproduces the result of~\cite{chen2008witten}.
\end{cor}
\section*{Acknowledgement}
This paper is a short version of the author's Master thesis. 
I am grateful to Doctor Qingtao Chen for a lot of discussion and help on the Witten genus. I also thanks Doctor Honglu Fan for inspiring suggestions on toric varieties.  I would like to thank Professor Ana da Silva for her kind help in clarifying some doubts about symplectic toric variety. Special thank goes to Professor Giovanni Felder. Thanks to his patient discussion and sharp questions, I managed to correct the false understandings of some black-boxes.

\end{document}